\documentclass{article}

\usepackage{amsmath,amsthm,enumerate,eucal,pgf,tikz}
\usepackage[utf8]{inputenc}

\DeclareMathOperator{\Q}{\mathbf{Q}}
\DeclareMathOperator{\F}{\mathbf{F}}
\DeclareMathOperator{\Z}{\mathbf{Z}}
\DeclareMathOperator{\C}{\mathbf{C}}
\newcommand\HH{\mathrm{H}}

  \renewenvironment{thebibliography}[1]{%
    \begin{oldthebibliography}{#1}%
      \setlength{\parskip}{0ex}%
      \setlength{\itemsep}{0.1ex}%
      \setlength{\labelwidth}{1.4cm} 
  }%
  {%
    \end{oldthebibliography}%
  }

\usepackage[T1]{fontenc}
\usepackage[scaled]{beramono,berasans}
\usepackage[charter]{mathdesign}
\usepackage{microtype}

\title{Computing overconvergent forms for small primes}
\author{Jan Vonk}
\date{}

\newtheorem{theorem}{Theorem}
\newtheorem{lemma}{Lemma}

\newtheorem{thmx}{Theorem}

\begin{document}

\maketitle
\begin{abstract}
In this note, we construct explicit bases for spaces of overconvergent $p$-adic modular forms when $p=2,3$ and study their stability under the Atkin operator. The resulting extension of the algorithms of Lauder \cite{Lau1,Lau2} is illustrated  with computations of slope sequences of some $2$-adic eigencurves and the construction of Chow-Heegner points on elliptic curves via special values of Rankin triple product L-functions.
\end{abstract}

\section*{Introduction}

\par Overconvergent modular forms were systematically developed by Katz in \cite{Kat1}, and have since played a prominent role in number theory. Many results in \cite{Kat1} require the existence of a lift of the Hasse invariant to $\Z_p$. In this note, we extend some of these results to situations where such a lift does not necessarily exist. 

\par The outline is as follows. In Section 1, we construct an explicit basis for spaces of overconvergent forms. As a result, we obtain a $\Z_p$-lattice which we compare to the one considered in \cite{Kat1}. In Section 2, we present a number of applications of these results by removing the restriction $p\geq 5$ from both Wan's quadratic bound for the Gouv\^ea-Mazur conjecture \cite{Wan1}, and Lauder's algorithms for computing characteristic series for $U_p$ and constructing Chow-Heegner points on elliptic curves \cite{Lau1,Lau2}. 

\section{Explicit bases and lattices}

\par In this section, we will construct explicit bases for $r$-overconvergent forms for any prime $p$, removing the restriction $p\geq 5$ from \cite[Proposition 2.6.2]{Kat1}. We then study how $U_p$ interacts with the lattice corresponding to the natural supremum-norm on the given basis. 

\par \textbf{Definitions. }Let $p$ be a prime and $K$ a finite extension of $\Q_p$ with valuation $v_p$, normalised such that $v_p(p)=1$, and valuation ring $R$. Set $(E,n)$ to denote $(E_{p-1},1)$ when $p\geq 5$, but $(E_4,4)$ when $p=2$ and $(E_6,3)$ when $p=3$. Here $E_k$ is the Eisenstein series of level $1$ and weight $k$. Note that $E$ is a lift to $\Z_p$ of the $n$-th power of the Hasse invariant, see Remark 2. Let $\mathcal{X} \rightarrow \mathrm{Spec}(\Z_p)$ be the compactified modular curve of level $\Gamma = \Gamma_1(N)$ for $N\geq 5$ coprime to $p$, with generic fibre $X$. Let $C$ be the closed subscheme of cusps on $\mathcal{X}$, and $\pi: \mathcal{U} \rightarrow \mathcal{X}$ the universal generalised elliptic curve with $\Gamma$-level structure. On $\mathcal{X}$, define the invertible sheaf $\omega := \pi_{*}\Omega^1_{\mathcal{U}/\mathcal{X}}(\log \pi^{-1}C)$, so that $E$ is a global section of $\omega^{\otimes n(p-1)}$. As in \cite[Section 1]{Col1}, we can define a rigid subspace $X_r$ of $X^{\mbox{\scriptsize{rig}}}$ for every $r \in \C_p$ whose points are exactly those $x$ such that $v_p(E_x) \leq v_p(r^n)$. We define the space of $r$-overconvergent modular forms of integer weight $k$ on $\Gamma$ to be $M^{\dagger,r}_k := \HH^0(X_r,\omega^{\otimes k})$. 

\par \textbf{Remark 1. }We can make similar definitions for arbitrary congruence subgroups $\Gamma$, and the computations in Section 2 are for $\Gamma_0(N)$. The justification lies in the fact that the arguments below can be done on the coarse moduli scheme $\mathcal{X}_0(N)$ over ${\Z_p}$ instead, working with the line bundles $\omega_k$ as in \cite[Lemma II.4.5]{Maz1}. A very careful analysis is given in \cite[Appendix]{BC05}.

\par \textbf{Remark 2. }We chose a particular value for $n$ for every prime, but everything continues to hold mutatis mutandis for $n$ an arbitrary power of $p$. For computational purposes there is an advantage in picking the smallest possible $n$ that assures a level $1$ lift, as we did above.

\subsection{Explicit bases}

\par We will now attempt to find an explicit basis for $M^{\dagger,r}_k$. As explained in \cite[Section 1]{Col1}, we can obtain $X_r$ as the Raynaud generic fibre of the completion along the special fibre of $Spec_X\left(\mbox{Sym}(\omega^{\otimes n(p-1)})/(E-r^n)\right)$. By reducing to the case of the open modular curve, which is affine, the analysis of the Leray spectral sequence in \cite[Theorem 2.5.1]{Kat1} shows that we can pull out the ideal and obtain
\begin{equation}\label{Quotient}
M^{\dagger,r}_k = \HH^0\left(X,\omega^{\otimes k} \otimes \mbox{Sym} (\omega^{\otimes n(p-1)})\right)/(E- r^n).
\end{equation}
Having this concrete description in hand, we now attempt to eliminate the relation $E = r^n$ by investigating the map given by multiplication by $E$ on modular forms as in \cite[Lemma 2.6.1]{Kat1}. The proof is nearly identical.

\begin{lemma}\label{Split}
Let $k \neq 1$, then the injection given by the multiplication by $E$-map
\[\HH^0\left(\mathcal{X},\omega^{\otimes k}\right) \stackrel{\times E}{\longrightarrow} \HH^0\left(\mathcal{X},\omega^{\otimes k + n(p-1)}\right)\]
splits as a map of $\Z_p$-modules. 
\begin{proof}The result is clear for $k\leq 0$. For $k\geq 2$, we have $\HH^1(\mathcal{X},\omega^{\otimes k}) = 0$ by computing the degree of $\omega$ as in \cite[Theorem 1.7.1]{Kat1}. We obtain the short exact sequence
\[0 \rightarrow \HH^0\left(\mathcal{X},\omega^{\otimes k}\right) \stackrel{\times E}{\longrightarrow} \HH^0\left(\mathcal{X},\omega^{\otimes k+n(p-1)}\right) \longrightarrow \HH^0\left(\mathcal{X},\mathcal{F}\right) \rightarrow 0,\]
where $\mathcal{F}$ is the quotient sheaf. This sequence remains exact after tensoring with $\F_p$. Indeed, $\mathcal{F}$ is flat over $\Z_p$ as $E$ is not identically $0$ in the special fibre, and since $\mathcal{F}$ is a skyscraper sheaf over $\F_p$ it follows that $\HH^1(\mathcal{X}_{\mathbf{F}_p},\mathcal{F}_{\mathbf{F}_p})=0$ and hence $\mathrm{Supp}\ \mathrm{R}^1f_{*}\mathcal{F}= \emptyset$, where $f: \mathcal{X}\rightarrow \mathrm{Spec}(\Z_p)$ is the defining morphism for $\mathcal{X}$. We conclude that $\HH^0\left(\mathcal{X},\mathcal{F}\right)$ is a free $\Z_p$-module, from which the conclusion follows.
\end{proof}
\end{lemma}

\par For every $i \geq 0$, use the above lemma to choose generators $\{a_{i,j} \}_{j}$ for a complement of the submodule $E\cdot \HH^0(\mathcal{X},\omega^{\otimes k + (i-1)n(p-1)})$ inside $\HH^0(\mathcal{X},\omega^{\otimes k + in(p-1)})$. This choice is non-canonical, but we will fix it once and for all in what follows. By running through the proof of \cite[Proposition 2.6.2]{Kat1}, one can check that the following theorem is a direct consequence of equation (\ref{Quotient}) and Lemma \ref{Split}.
\begin{thmx}\label{Basis}
The set $\left\{r^{ni}a_{i,j}E^{-i}\right\}_{i,j}$ is a basis for the $p$-adic Banach space $M^{\dagger,r}_k$.
\end{thmx}

\par \textbf{Remark 3. }Note that we have avoided the case $k = 1$, as the standard base change results are known to fail for many levels. However, in all our applications we can get the required result in weight $1$ by a simple application of Frobenius linearity of $U_p$ in the sense of \cite[Eqn. (3.3)]{Col2}, hence reducing the question to one in higher weight for which the results above hold. See also \cite[Section 2.2]{Lau1}.

\subsection{Comparing lattices}

We now have two integral structures on $M^{\dagger,r}_k$. The first, which we will call $\mathcal{B}_H(r)$, is the one defined by Katz \cite[Section 2.2]{Kat1} using the interpretation of modular forms as certain rules on test objects. It has the advantage of being well suited for analysing its interaction with various operators. It has the disadvantage of only being explicit (in the sense of Theorem A) when a lift of the Hasse invariant to $\Z_p$ exists. The second, which we will call $\mathcal{B}_E(r)$, is simply the collection of forms that have integral coordinates with respect to the chosen basis from Theorem A. It has the advantage of being explicit and computational, but the disadvantage of being non-canonical and hence having a rather mysterious interaction with $U_p$. We will now attempt to compare $\mathcal{B}_H$ and $\mathcal{B}_E$, in order to get the best of both worlds. 

\begin{lemma}\label{Lattices}
Assume that there exists a lift $H$ of the Hasse invariant to $\Z_p$, then we have
\[r^{n-1}\mathcal{B}_H(r) \subseteq \mathcal{B}_E(r).\]
\begin{proof}Fix a choice of complementary subspaces for both $E$ and $H$ as above. Let $f \in \mathcal{B}_H(r)$, then by Theorem A we can write $f = \sum_{i\geq 0}\ \  r^{i}a_iH^{-i}$, where $a_i$ is in the $i$-th complementary subspace for $H$. We rewrite this as 
\begin{equation}\label{HForm}f =  \sum_{i \geq 1} \sum_{j=0}^{n-1} r^{ni-j}\ a_{ni-j}\  H^{-(ni-j)} = \sum_{i \geq 1}\ r^{ni}H^{-ni} \left(\sum_{j=0}^{n-1} \ \ \ r^{-j}\ a_{ni-j}\  H^{j}\right). \end{equation}
The inner sum in the above expression is guaranteed to be in $\HH^0(\mathcal{X},\omega^{\otimes k+ni})$ when multiplied by $r^{n-1}$. We can decompose this multiple as $\sum_{m = 0}^{i}b_mE^m$, where $b_m$ is in the $m$-th complementary subspace for $E$. Recall that $n$ is a power of $p$, from which we get $E \equiv H^{n} \pmod{pn}$. If we substitute all this into (\ref{HForm}), we obtain that $r^{n-1}f \in \mathcal{B}_E(r)$ as desired.
\end{proof}
\end{lemma}

With the aid of this lemma, we now investigate the interaction of our explicit lattice $\mathcal{B}_E$ with the operators $U_p$ and multiplication by $G := \frac{E}{V_pE}$, where $V_p$ is the Frobenius operator defined in \cite[Section 2]{Col1}. Both operators will play a crucial role in the applications.

\begin{thmx}\label{Stable}
Let $v_p(r) < \frac{1}{p+1}$, then we have 
\[r^{n-1}pU_p\left( \mathcal{B}_E(r)\right) \subseteq \mathcal{B}_E(r^p) \hspace{.4cm} \mbox{ and } \hspace{.4cm} G\cdot \mathcal{B}_E(r) \subseteq \mathcal{B}_E(r).\] 
\begin{proof}Assume first that there exists a lift of the Hasse invariant to $\Z_p$. The first statement follows immediately from Lemma \ref{Lattices} and the inclusion $pU_p\left(\mathcal{B}_H(r)\right) \subseteq \mathcal{B}_H(r^p)$, which is \cite[Proposition II.3.6]{Gou1}. For the second statement, we will use that $G \equiv 1 \pmod {pr^{n-p-2}}$ when $v_p(r) < 1/(p+1)$. For $p \geq 5$, this is exactly \cite[Lemma 2.1]{Wan1}. For $p\leq 3$, we can check this directly from the formulae in \cite[Section 3]{Cal1}. It follows that $G \cdot \mathcal{B}_E \subseteq \mathcal{B}_E + r^{n-1}\mathcal{B}_H \subseteq \mathcal{B}_E$.

\par If no lift of the Hasse invariant to $\Z_p$ exists, add two additional level structures that both assure existence and intersect trivially, see \cite[Section 6]{Col1} and \cite[Section B2]{Col2}. The result now follows from taking intersections on the level of Katz expansions.
\end{proof}
\end{thmx}

\section{Applications}

We now sketch how Theorems A and B enable us to generalise previous work in the literature due to Wan \cite{Wan1}, Lauder \cite{Lau1} and Darmon-Lauder-Rotger \cite{DLR1}. We work with $\Gamma = \Gamma_0(N)$ for computational simplicity when appropriate, see Remark 1.

\subsection{The Gouv\^ea-Mazur conjecture}

\par An enormous amount of arithmetic information is encoded in the \textit{slopes} of overconvergent modular forms, which are the valuations of their $U_p$-eigenvalues. One of the consequences of the theory of Coleman \cite{Col2} is that for any $\alpha > 0$, there exists a smallest integer $N_{\alpha}$ with the following property: If $k_1,k_2 \in \Z$ such that $k_1 \equiv k_2 \mod p^{N_{\alpha}}(p-1)$, then the collection of slopes $\leq \alpha$ in weights $k_1$ and $k_2$ agree, with multiplicities. Gouv\^ea and Mazur conjectured in \cite{GM1} that $N_{\alpha} \leq \lfloor \alpha \rfloor$, to which a counterexample was given in \cite{BC1}. However, Wan \cite{Wan1} exhibits an explicit quadratic upper bound for $N_{\alpha}$, provided that $p\geq 5$. We will now remove this restriction on $p$. 

\par Wan's analysis relies on a good knowledge of an explicit basis, along with an understanding of how $U_p$ and $G$ act on the integral lattice. This is exactly the content of Theorems A and B, making the proof a straightforward adaptation of the methods in \cite{Wan1}. We estimate the size of the coefficients of the characteristic series $P_k(t)$ of $U_p$ on the space $M_k^{\dagger,r}$. This is done by analysing the entries of the matrix of $U_p$ on the Katz basis. After twisting $U_p$ by $E$, Theorem B enables us to do this uniformly with respect to variations of the weight. 

\par \textbf{Notation. }Choose generators $a_{u,v}$ for the $u$-th complementary subspace, giving rise to a Katz basis $e_{u,v} = r^{nu}a_{u,v}E^{-u}$. Multiplication by $E^j$ defines an isomorphism $M_k^{\dagger,r} \rightarrow M_{k+jn(p-1)}^{\dagger,r}$, so we conclude by an application of Frobenius linearity of $U_p$ \cite[Eqn. (3.3)]{Col1} that
$P_{k+jn(p-1)}(t)$ equals the characteristic series of $U_p \circ G^j$ on $ M_k^{\dagger,r}$, where we recall that $G = E(V_pE)^{-1}$. We write 
\[U \circ G^j(e_{u,v}) = \sum_{w,z} A_{u,v}^{w,z}(j)\ e_{w,z},\]
for some $A_{u,v}^{w,z} \in K$. The following lemma estimates the size of these numbers, independently of $j$.

\begin{lemma}[Cfr. \protect{\cite[Lemma 3.1]{Wan1}}]\label{Matrix}
We have \[v_p\left(A_{u,v}^{w,z}(j) \right) \geq wn(p-1)v_p(r) - 1 - v_p(r)(n-1).\]
\begin{proof}It follows from Theorem B that 
\[U \circ G^j (e_{u,v}) = \frac{1}{r^{n-1}p} \sum_w \frac{r^{npw}}{E^w}b_w(u,v,j) = \frac{1}{r^{n-1}p} \sum_w r^{n(p-1)w} \frac{r^{nw}}{E^w}b_w(u,v,j),\]
where $b_w(u,v,j)$ is in the $j$-th complementary subspace, and hence an integral combination of the $a_{w,z}$. This gives us the desired bound on $A_{u,v}^{w,z}(j)$.
\end{proof}
\end{lemma}

\par The key observation is that the above lower bound is independent of $j$. After taking determinants, we obtain a lower bound on the coefficients of $P_{k+jn(p-1)}(t)$, again independent of $j$. Wan now proceeds by proving a very general reciprocity lemma on Newton polygons, which allows him to transform the lower bound for $P_k(t)$ into an upper bound for $N_{\alpha}$. The analysis goes through without modifications, and using Wan's results we deduce from Lemma \ref{Matrix} that
\begin{theorem}
There is an explicitly computable quadratic polynomial $P \in \mathbf{Q}[x]$ such that $N_{\alpha} \leq P(\alpha)$.
\end{theorem}

\subsection{Buzzard's slope conjectures}

\par Using the explicit bases for $M^{\dagger,r}_k$ in \cite{Kat1} and the $p$-adic estimates in \cite[Lemma 3.1]{Wan1}, Lauder presents an algorithm to compute the characteristic series of $U_p$ when $p\geq 5$ in \cite{Lau1}. By an application of Coleman's trick \cite[Eqn. (3.3)]{Col1}, it is particularly useful when $k$ becomes very large. Given the theory above, it is straightforward to remove the restriction on $p$. The code for our extension to small primes can be found on the author's webpage. In what follows, we will explicitly compute some examples. All computations were a matter of seconds on a standard laptop. 

\par In \cite{Buz1}, Buzzard made very precise conjectures on the sequence of slopes for $M^{\dagger,r}_k$ on $\Gamma_0(N)$, and gives a precise conjectural recipe when $p$ is $\Gamma_0(N)$-regular. This is a condition which essentially ensures that the slopes at small weights are as small as Hida theory allows them to be. For a precise definition and a reformulation in terms of Galois representations, see \cite[Section 1]{Buz1}.

\par \textbf{Example 1. }We compute that the first few slopes of $U_3$ acting on $M^{\dagger}_{278}\left(\Gamma_0(41)\right)$ are 
\[\mathbf{0}_{12},\mathbf{1}_{14},\mathbf{3}_{48}, \mathbf{6}_{14}, \mathbf{7}_{22}, \mathbf{8}_{6}, \mathbf{9}_{22}, \mathbf{10}_{14}, \mathbf{12}_{48}, \mathbf{14}_{14}, \mathbf{16}_{22}, \mathbf{17}_6, \mathbf{18}_{22},\ldots \]
where the subscripts denote multiplicities. We check that $3$ is $\Gamma_0(41)$-regular, and that the slopes agree with Buzzard's prediction. Note that this slope sequence equals the one in weight $8$ for all the terms we display here, suggesting a very strong form of the Gouv\^ea-Mazur conjecture. 

\par \textbf{Example 2. }To illustrate a case where regularity fails in a striking way, we compute the first few slopes of $U_2$ acting on $M^{\dagger}_{10}\left(\Gamma_0(89)\right)$, where as before the subscripts denote multiplicities:
\[\mathbf{0}_{16}, \mathbf{1}_{22}, \mathbf{2}_{22}, \mathbf{14/5}_5, \mathbf{3}_1, \mathbf{4}_{68}, \mathbf{9/2}_4, \mathbf{6}_1, \mathbf{31/5}_5, \mathbf{7}_{22}, \mathbf{8}_{22}, \mathbf{9}_{30}, \mathbf{10}_{22}, \mathbf{21/2}_{16}, \mathbf{12}_{52}, \ldots\]
The appearance of denominators as large as $5$ does not seem to have been recorded before. Note that by Coleman \cite{Col1}, the overconvergent forms giving rise to these denominators are in fact even classical. 

\par \textbf{Example 3. }A more systematic computation of $2$-adic overconvergent forms of levels $\Gamma_0(53)$ and $\Gamma_0(61)$ suggests a remarkable relationship between the corresponding eigencurves, for which we have no explanation. The table below lists the first few entries of the $2$-adic slope sequences in weights $14$ and $16$. 

\begin{table}[!ht]
\centering
\begin{tabular}{c|l}
 & $k = 14$ \\
\hline
$\Gamma_0(53)$ & $\mathbf{0}_{10}, \mathbf{1}_{13}, \mathbf{2}_{23}, \mathbf{4}_{13}, \mathbf{6}_{59}, \mathbf{9}_{13}, \mathbf{11}_{23}, \mathbf{12}_{13}, \mathbf{13}_{18}, \mathbf{14}_{13}, \mathbf{29/2}_{10}, \mathbf{16}_{18}, \mathbf{17}_{13}, \mathbf{18}_{23}, \mathbf{21}_{13}, \ldots$ \\
$\Gamma_0(61)$ & $\mathbf{0}_{12}, \mathbf{1}_{15}, \mathbf{2}_{25}, \mathbf{4}_{15}, \mathbf{6}_{69}, \mathbf{9}_{15}, \mathbf{11}_{25}, \mathbf{12}_{15}, \mathbf{13}_{22}, \mathbf{14}_{15}, \mathbf{29/2}_{10}, \mathbf{16}_{22}, \mathbf{17}_{15}, \mathbf{18}_{25}, \mathbf{21}_{15}, \ldots$\\
\vspace{.2cm} & \\
 & $k = 16$ \\
 \hline
$\Gamma_0(53)$ & $\mathbf{0}_{10}, \mathbf{1}_{13}, \mathbf{3/2}_{10}, \mathbf{3}_{31}, \mathbf{17/3}_3, \mathbf{6}_1, \mathbf{7}_{67}, \mathbf{15/2}_2, \mathbf{9}_1, \mathbf{28/3}_3, \mathbf{12}_{31}, \mathbf{27/2}_{10}, \mathbf{14}_{13}, \mathbf{15}_{18}, \mathbf{16}_{13}, \ldots$\\
$\Gamma_0(61)$ & $\mathbf{0}_{12}, \mathbf{1}_{15}, \mathbf{3/2}_{10}, \mathbf{3}_{37}, \mathbf{17/3}_3, \mathbf{6}_1, \mathbf{7}_{78},\hspace{.28cm}\mathbf{8}_1, \hspace{.28cm}\mathbf{9}_1, \mathbf{28/3}_{3}, \mathbf{12}_{37}, \mathbf{27/2}_{10}, \mathbf{14}_{15}, \mathbf{15}_{22}, \mathbf{16}_{15}, \ldots$\\
\end{tabular}
\end{table}

This computation was carried out to a large precision, and for a much larger range of weights. We chose to include the start of the sequence for $k=14,16$ as it illustrates the general behaviour rather well. The set of slopes, without multiplicities, seems to agree for both levels in all weights, with the exception of a small deviation. This deviation, if it occurs, seems to come from the $2$-stabilisations of the largest classical cuspidal slope of level $N$. 

\subsection{Chow-Heegner points on elliptic curves}

\par We will now use our results to explicitly perform some Heegner-type point constructions on elliptic curves, following the theory in \cite{DRI} and the algorithm of \cite{Lau2}, which was conditional on $p \geq 5$.

\par Let $p$ be a prime and $E/\Q$ an elliptic curve of conductor $N$, associated to the $p$-ordinary form $f \in S_2^{new}(\Gamma_0(N))$, and let $g$ be any other weight $2$ newform which is $p$-ordinary. As explained in \cite[Section 1]{Lau2}, we can deduce from \cite[Theorem 1.3]{DRI} that there exists a point $P_g \in E(\Q)$ such that
\begin{equation}\label{DR} \log(P_g) = 2d_g\ \cdot\  \frac{\mathcal{E}_0(g) \ \mathcal{E}_1(g)}{\mathcal{E}(g,f,g)}\ \cdot\ \mathcal{L}_p (\mathbf{g},\mathbf{f}, \mathbf{g})(2,2,2)\ ,
\end{equation}
where $\log$ denotes the formal $p$-adic logarithm on $E$, $d_g$ is an integer described in \cite[Remark 3.1.3]{DDLR}, the $\mathcal{E}$-factors are computable quadratic numbers depending only on the $p$-th coefficients of $f$ and $g$ described in \cite[Theorem 1.3]{DRI}, and $\mathcal{L}_p(\mathbf{g},\mathbf{f}, \mathbf{g})$ is the Rankin triple product $p$-adic L-function associated to the Hida families $\mathbf{f},\mathbf{g}$ through $f,g$ respectively. 

\par The crux in computing the special value of the Rankin triple product $p$-adic L-function is the efficient computation of the $U_p$-operator on the space $M^{\dagger,r}_k$. The previous subsection removed the restriction $p\geq 5$ from the algorithm in \cite{Lau1} to compute this action, and it is now straightforward to compute the desired special value of the Rankin triple product $p$-adic L-function, as described in detail in \cite{Lau2}. We have implemented a version in \texttt{Magma} that works for all $p$, which is available on the author's webpage. Let us turn to some numerical examples.

\par \textbf{Example 1. }Let $E : y^2 + xy = x^3 - x^2 - x + 1$ be the rank $1$ elliptic curve of conductor $58$, with associated newform $f$, and let $g$ be the unique newform on $\Gamma_0(58)$ different from $f$. Both $f$ and $g$ are $2$-ordinary. Letting $P = (0,1)$ be a generator for $E(\Q)$, we compute that
\[\mathcal{L}_2(\mathbf{g},\mathbf{f}, \mathbf{g})(2,2,2) \equiv 3\log_E(P) \pmod {2^{200}},\]
as predicted by the theory in \cite{DRI}.  

\par \textbf{Remark 5. }When the Tate module of $E_{\Q}$ is wildly ramified at $2$ or $3$, we might wonder whether the Chow-Heegner point construction just described continues to work. The associated newform $f$ will be of infinite slope, so we lack a notion of Hida or Coleman family passing through $f$. It is therefore not obvious whether the theoretical framework of \cite{DRI} will generalise to such a setting. Nonetheless, we are often able to run our extension of Lauder's algorithms \cite{Lau1} and \cite{Lau2}, and recover a rational point on $E$, as the following examples show.

\par \textbf{Example 2a. }Let $E : y^2 + y = x^3 + 9x - 10$, which is an elliptic curve over $\Q$ of conductor $4617=3^5\cdot 19$ and rank $1$. Let $f$ be the associated newform. Let $g = q - 2q^3 - 2q^4 +3q^5-q^7 + \ldots$ be the unique cuspidal newform of weight $2$ on $\Gamma_0(19)$. Despite $f$ being of infinite $3$-adic slope, we can run the computation and find a numerical value for $\mathcal{L}_2(\mathbf{g},``f", \mathbf{g})(2,2,2)$. We find that 
\[\mathcal{L}_3(\mathbf{g},``f", \mathbf{g})(2,2,2) \equiv t\cdot \log_E(P) \pmod{3^{200}}\ \ \mbox{   where   } \ \ 2t^2 + 48t + 729=0,\]
where $P = (4,9)$ is a generator of $E(\Q)$. The fact that both quantities are related by a quadratic number $t$ of small height suggests that a more general analogue of the theory for ordinary forms in \cite{DRI}, and more specifically equation (\ref{DR}), might exist.

\par \textbf{Example 2b. }Let $E : y^2 = x^3 + x^2 - 62893x - 6091893$, which is an elliptic curve over $\Q$ of rank $1$ and conductor $15104 = 2^8\cdot 59$. Let $f$ be its associated newform, and let $g = q - q^2 - q^3 + q^4 - 3q^5 + \ldots$ be the newform of level $118$ associated to the elliptic curve with Cremona label \texttt{118.a1}. Note that $g$ is $2$-ordinary. We compute that 
\[\mathcal{L}_2(\mathbf{g},``f", \mathbf{g})(2,2,2) \equiv 6\log_E(P) \pmod {2^{100}},\]
where $P = (20821, 3004216)$ is a generator of $E(\Q)$. As in the previous example, this suggests that an analogue of equation (\ref{DR}) holds for $f$ of infinite slope. Note that this would work the other way: once we compute the value of $\mathcal{L}_2(\mathbf{g},``f", \mathbf{g})(2,2,2)$, we can use a formal exponentiation routine as in \cite{Lau2} to recover a point of infinite order, which is of considerable height in this example.

\par \textbf{Remark 6. }Chow-Heegner points have a well-understood geometric origin and can also be constructed by complex analytic methods, see \cite{DRS} and \cite{DDLR}. For an application of triple product $p$-adic L-functions, and the methods in this paper, to the $p$-adic construction of points in more mysterious settings we refer the reader to \cite{DLR1}.

\section*{Acknowledgements}

\par We wish to thank Alan Lauder for suggesting this problem and providing generous assistance with the computations, Kevin Buzzard and David Loeffler for their useful comments and suggestions. The author is supported by the EPSRC and St. Catherine's College, Oxford. 

\bibliographystyle{alpha}

\end{document}